\documentclass[11pt]{amsart}
\usepackage[mathscr]{eucal}
\usepackage{amsfonts}
\usepackage{amsmath}
\usepackage{amsthm}
\usepackage{amssymb}
\usepackage{amscd}
\usepackage{latexsym}
\usepackage{slashbox}
\usepackage{lscape}

\theoremstyle{plain}
\newtheorem{definition}[equation]{Definition}

\newtheorem{lemma}[equation]{Lemma}
\newtheorem{proposition}[equation]{Proposition}
\newtheorem{theorem}[equation]{Theorem}

\theoremstyle{remark}
\newtheorem{remark}[equation]{Remark}
\newtheorem{notation}[equation]{Notation}

\numberwithin{equation}{subsection}
\renewcommand{\mathfrak}{\mathcal}

\begin{document}
\title{The calculus structure of the Hochschild homology/cohomology of preprojective algebras of Dynkin quivers }
\author{Ching-Hwa Eu}
\address{Department of Mathematics, Massachusetts Institute of
Technology, Cambridge, MA 02139, U.S.A.}
\email{ceu@math.mit.edu}
\maketitle
\pagestyle{myheadings}
\markboth{Ching-Hwa Eu}{The calculus of the Hochschild (co)homology of preprojective algebras of Dynkin quivers}

\section{Introduction}
The Hochschild cohomology ring of any associative algebra, together with the Hochschild homology, forms a structure of calculus. This was proved in \cite{DGT}. In this paper, we compute the calculus structure for the preprojective algebras of Dynkin quivers over a field of characteristic zero, using the Batalin-Vilkovisky structure of the Hochschild cohomology. Together with the results of \cite{CBEG}, where the Batalin-Vilkovisky structure is computed for non-ADE quivers (and the calculus can be easily computed from that), this work gives us a complete description of the calculus for any quiver.

The Hochschild homology and cohomology spaces and the duality between them was established in \cite{EE2}. The cup product structure of the Hochschild cohomology was computed in \cite{ES2} for the quivers of type $A$ and \cite{Eu2} for the ones of type $D$ and $E$. We use the notations from these papers and give the calculus structure in terms of the bases which were defined there.

First, we compute the Connes differential on Hochschild homology by using the Cartan identity. Since it turns out this differential makes the Hochschild cohomology ring a Batalin-Vilkovisky-algebra, this gives us an easy way to compute the Gerstenhaber bracket and the contraction map. Then we use the Cartan identity to compute the Lie derivative.

{\bf{Acknowledgements.}} C. Eu wants to thank his advisor P. Etingof for useful discussions and V. Dolgushev for his explainations about calculus.

\section{Preliminaries}

\subsection{Quivers and path algebras}
Let $Q$ be a quiver of ADE type with vertex set $I$ and $|I|=r$. 
We write $a\in Q$ to say that $a$ is an arrow in $Q$. 

We define $Q^*$ to be the quiver obtained from $Q$ by reversing
all of its arrows. We call $\bar Q=Q\cup Q^*$ the \emph{double}
of $Q$. 

Let $C$ be the adjacency matrix corresponding to the
quiver $\bar Q$. 

The concatenation of arrows generate the \emph{nontrivial
paths} inside the quiver $\bar Q$. We define $e_i$, $i\in I$ to
be the \emph{trivial path} which starts and ends at $i$. The
\emph{path algebra} $P_{\bar Q}=\mathbb{C}\bar Q$ of $\bar Q$ over $\mathbb{C}$ is the $\mathbb{C}$-algebra with basis the paths in $\bar Q$ and the product $xy$ of two paths $x$ and $y$ to be their concatenation if they are compatible and $0$ if not. We define the \emph{Lie bracket} $[x,y]=xy-yx$.

Let $R=\oplus_{i\in I}\mathbb{C}e_i$. Then $R$ is a commutative semisimple
algebra, and $P_Q$ is naturally an
$R$-bimodule.
\begin{notation}
 Throughout this paper, we write ''$\otimes=\otimes_R$''.
\end{notation}

\subsection{Frobenius algebras}
Let $A$ be a finite dimensional unital $\mathbb{C}-$algebra. We
call it Frobenius if there is a linear function $f:A\rightarrow\mathbb{C}$, such that the form $(x,y):=f(xy)$ is nondegenerate, or, equivalently, if there exists an isomorphism $\phi:A\stackrel{\simeq}{\rightarrow}A^*$ of left $A-$modules: given $f$, we can define $\phi(a)(b)=f(ba)$, and given $\phi$, we define $f=\phi(1)$.

If $\tilde f$ is another linear function satisfying the same
properties as $f$ from above, then $\tilde f(x)=f(xa)$ for some
invertible $a\in A$. Indeed, we define the form $\{a,b\}=\tilde f(ab)$. Then $\{-,1\}\in A^*$, so there is an $a\in A$, such that $\phi(a)=\{-,1\}$. Then $\tilde f(x)=\{x,1\}=\phi(a)(x)=f(xa)$.

\subsection{The Nakayama automorphism}
Given a Frobenius algebra $A$ (with a function $f$ inducing a
bilinear form $(-,-)$ from above), the automorphism
$\eta:A\rightarrow A$ defined by the equation $(x,y)=(y,\eta(x))$
is called the \emph{Nakayama automorphism} (corresponding to
$f$). 

\subsection{The preprojective algebra}
Given a quiver $Q$, we define the \emph{preprojective
algebra} $\Pi_Q$ to be the quotient of the path algebra $P_{\bar Q}$ by
the relation $\sum\limits_{a\in Q}[a,a^*]=0$. 

Given a path $x$, we write $x^*$ for the path obtained from $x$ by reversing all arrows.

From now on, we write $A=\Pi_{Q}$.

\subsection{Graded spaces}

Let $M=\oplus_{d\geq0}M(d)$ be a $\mathbb Z_+$-graded
vector space, with finite dimensional homogeneous subspaces. 
We denote by $M[n]$ the same space with
grading shifted by $n$. The graded dual space $M^*$ is defined by the
formula $M^*(n)=M(-n)^*$.

\subsection{Root system parameters}

Let $w_0$ be the longest element of the Weyl group $W$ of $Q$. 
Then we define $\nu$ to be the involution of $I$, such
that $w_0(\alpha_i)=-\alpha_{\nu(i)}$ (where $\alpha_i$ is the
simple root corresponding to $i\in I$). It turns out that
$\eta(e_i)=e_{\nu(i)}$ (\cite{S}; see \cite{ES2}).

Let $m_i$, $i=1,...,r$, be the exponents of the root system
attached to $Q$, enumerated in the increasing order. 
Let $h=m_r+1$ be the Coxeter number of $Q$. 

Let $P$ be the permutation matrix corresponding to the involution
$\nu$. Let $r_+=\dim\ker(P-1)$ and $r_-=\dim\ker(P+1)$.
Thus, $r_-$ is half the number of vertices which are not fixed by
$\nu$, and $r_+=r-r_-$.

$A$ is finite dimensional, and the following Hilbert series is known from \cite[Theorem 2.3.]{MOV}:

\begin{equation}\label{Hilbert series matrix}
H_A(t)=(1+Pt^h)(1-Ct+t^2)^{-1}.
\end{equation}

We see that the top degree of $A$ is $h-2$, and for the top degree $A^{top}$ part we get the following decomposition in $1$-dimensional submodules:

\begin{equation}\label{topdegree}
A^{top}=A(h-2)=\bigoplus_{i\in I}e_iA(h-2)e_{\nu(i)}
\end{equation}

\section{Hochschild cohomology and homology}
The Hochschild cohomology and homology spaces of $A$ were computed in \cite{EE2}. We recall the results:
\begin{definition}
We define the spaces 
\begin{eqnarray*}
U&=&\oplus_{d<h-2}HH^0(A)(d)[2],\\
L&=&HH^0(A)(h-2),\\
K&=&HH^2(A)[2],\\
Y&=&HH^6(A)(-h-2).
\end{eqnarray*}
\end{definition}

\begin{theorem}\label{t1}
 For the Hochschild cohomology spaces, we have the following natural isomorphisms:
\begin{align*}
{HH^0}(A)&=U[-2]\oplus L[h-2],\\
{HH^1}(A)&=U[-2],\\
{HH^2}(A)&=K[-2],\\
{HH^3}(A)&=K^*[-2],\\
{HH^4}(A)&=U^*[-2],\\
{HH^5}(A)&=U^*[-2]\oplus Y^*[-h-2],\\
{HH^6}(A)&=U[-2h-2]\oplus Y[-h-2],
\end{align*}
and ${HH^{6n+i}}(A)={HH^i}(A)[-2nh]\,\forall
i\geq1$.
\end{theorem}

\begin{theorem}\label{t2}
The Hochschild homology spaces of $A$, as graded spaces, are
as follows:
\begin{align*}
{HH_0}(A)=R,\\
{HH_1}(A)=U,\\
{HH_2}(A)=U\oplus Y[h],\\
{HH_3}(A)=U^*[2h]\oplus Y^*[h],\\
{HH_4}(A)=U^*[2h].\\
{HH_5}(A)=K[2h],\\
{HH_6}(A)=K[2h],\\
\end{align*}
and ${HH_{6n+i}}(A)={HH_i}(A)[2nh]\,\forall i\geq1$. 
\end{theorem}
In \cite{EE2}, an isomorphism $HH_\bullet(A)=HH^{8-\bullet}(A)[2h+2]$ was introduced. However, because of the periodicity of the Schofield resolution (with period $6$), we get for every $m\geq0$ an isomorphism 
\begin{equation}\label{H-dual}
\mathbb{D}:HH_\bullet(A)\stackrel{\sim}{\rightarrow}HH^{6m+2-\bullet}(A)[2mh+2]
\end{equation}

\begin{notation}\label{bases}
 In and \cite{Eu2}, the basis elements $z_k\in U[-2]\subset HH^0(A)$, $\omega_k\in L[h-2]$, $\theta_k\in U[-2]=HH^1(A)$, $f_k\in K[-2]=HH^2(A)$, $h_k\in K^*[-2]$, $\zeta_k\in U^*[-2]=HH^4(A)$, $\psi_k\in U^*[-2]\subset HH^5(A)$, and $\varepsilon_k\in Y^*[-h-2]$ were introduced. 
 
Our notations are taken from \cite{Eu2} and are different from those in \cite{ES2}:\\
$z_{k}$ in \cite{ES2} corresponds to $z_{2k}$ in our notations, \\
$z_m$ corresponds to $\omega_{\frac{h-3}{2}}$,\\
$g_k$ corresponds to $\theta_{2k}$,\\
$\psi_k$ corresponds to $\zeta_{l-k}$ and $\zeta_k$ corresponds to $\theta_{l-k}$ for $l=\left\{\begin{array}{cc}h-2&Q=A_n,\,n\,odd,\\h-3&Q=A_n,\,n\,even\end{array}\right.$
 For $c_k\in HH^i(A)$, $0\leq i\leq 5$, we write $c_k^{(s)}$ for the corresponding cocycle in $HH^{i+6s}$. 
We write $c_{k,t}$ for a cycle in $HH_{j+6t}$, $0\leq j\leq5$ which equals $\mathbb{D}^{-1}(c_k^{(s)})$.\\

We also introduced the maps $\alpha:K\rightarrow K^*$ and $\beta:Y^*\rightarrow Y$ there.
\end{notation}

\section{The calculus structure of the preprojective algebra}
We recall the definition of the calculus.
\subsection{Definition of calculus}
\begin{definition}\emph{(Gerstenhaber algebra)} 
A graded vector space $\mathcal{V}^\bullet$ is a Gerstenhaber algebra if it is equipped with a graded commutative and associative product $\wedge$ of degree 0 and a graded Lie bracket $[,]$ of degree $-1$. These operations have to be compatible in the sense of the following Leibniz rule

\begin{equation}
 [\gamma,\gamma_1\wedge\gamma_2]=[\gamma,\gamma_1]\wedge\gamma_2+(-1)^{k_1(k+1)}\gamma_1\wedge[\gamma,\gamma_2],
\end{equation}
where $\gamma\in\mathcal{V}^k$ and $\gamma_1\in\mathcal{V}^{k_1}$.
\end{definition}
We recall from \cite{CST} that
\begin{definition}\emph{(Precalculus)} 
A precalculus is a pair of a Gerstenhaber algebra $(\mathcal{V}^\bullet,\wedge,[,])$ and a graded vector space $\mathcal{W}^\bullet$ together with 
 \begin{itemize}
  \item a module structure $\iota_\bullet:\mathcal{V}^\bullet\otimes\mathcal{W}^{-\bullet}\rightarrow\mathcal{W}^{-\bullet}$ of the graded commutative algebra $\mathcal{V}^\bullet$ on $\mathcal{W}^{-\bullet}$.
  \item an action $\mathcal{L_\bullet}:\mathcal{V}^{\bullet+1}\otimes\mathcal{W}^{-\bullet}\rightarrow\mathcal{W}^{-\bullet}$ of the graded Lie algebra $\mathcal{V}^{\bullet+1}$ on $\mathcal{W}^{-\bullet}$ which are compatible in the sense of the following equations 
 \begin{equation}
  \iota_a\mathcal{L}_b-(-1)^{|a|(|b|+1)}\mathcal{L}_b\iota_a=\iota_{[a,b]},
 \end{equation}
and 
\begin{equation}
 \mathcal{L}_{a\wedge b}=\mathcal{L}_a\iota_b+(-1)^{|a|}\iota_a\mathcal{L}_b.
\end{equation}
 \end{itemize}
\end{definition}
\begin{definition}\emph{(Calculus)} 
 A calculus is a precalculus $(\mathcal{V}^\bullet,\mathcal{W}^\bullet,[,],\wedge,\iota_\bullet,\mathcal{L}_\bullet)$ with a degree 1 differential $d$ on $\mathcal{W}^\bullet$ such that the \emph{Cartan identity},
 \begin{equation}\label{Cartan}
  \mathcal{L}_a=d\iota_a-(-1)^{|a|}\iota_ad,
 \end{equation}
 holds.
\end{definition}
Let $A$ be an associative algebra. The contraction of the Hochschild cochain $P\in C^k(A,A)$ with the Hochschild chain $(a_0,a_1,\ldots,a_n)$ is defined by
 
\begin{equation}
 I_P(a_0,a_1,\ldots,a_n)=
 \left\{
 \begin{array}{cc}
(a_0P(a_1,\ldots,a_k),a_{k+1},\ldots,a_n) & n\geq k,\\
0   &   \mbox{else}.
\end{array}
\right.
\end{equation}

We have
\begin{proposition}\emph{\textbf{(Yu. Daletski, I. Gelfand and B. Tsygan \cite{DGT})}}
 The contraction $I_P$ together with the Connes differential, the Gerstenhaber bracket, the cup product and the action of cochains on chains (\cite[(3.5), page 46]{D}) induce on the pair $(HH^\bullet(A,A), HH_\bullet(A,A))$ a structure of calculus.
\end{proposition}

\section{Results about the calculus structure of the Hochschild cohomology/homology of preprojective algebrasof Dynkin quivers}

We state the results in terms of the bases of $HH^\bullet(A)$ and $HH_\bullet(A)$ which were introduced in Notation \ref{bases}:

\begin{theorem}
The calculus structure is given by tables 1, 2, 3 and the \emph{Connes differential $B$}, given as follows

 The Connes differential $B$ is given as follows:
 \begin{eqnarray*}
 B_{1+6s}(\theta_{k,s})&=&(1+\frac{k}{2}+sh)z_{k,s},\\
 B_{2+6s}(\omega_{k,s})&=&(\frac{1}{2}+s)h\beta^{-1}(\omega_{k,s}),\\
 B_{2+6s}(z_{k,s})&=&0,\\
 B_{3+6s}(\psi_{k,s})&=&((s+1)h-1-\frac{k}{2})\zeta_{k,s},\\
 B_{3+6s}(\varepsilon_{k,s})&=&0,\\
 B_{4+6s}&=&0,\\
 B_{5+6s}(h_{k,s})&=&(s+1)h\alpha^{-1}(h_{k,s}),\\
 B_{6+6s}&=&0.
 \end{eqnarray*}
\end{theorem}

\begin{landscape}
\thispagestyle{empty}
\begin{table}
\begin{tabular}[b]{|c||c|c|c|c|c|c|c|c|}
\hline
\backslashbox{$a$}{$b$}&$\theta_{l,t}$&$\omega_{l,t}$&$z_{l,t}$&$\psi_{l,t}$&$\varepsilon_{l,t}$&$\zeta_{l,t}$&$h_{l,t}$&$f_{l,t}$\\\hline\hline
$z_k^{(s)}$&$(z_k\theta_l)_{t-s}$&$\delta_{k0}\omega_{l,t-s}$&$(z_kz_l)_{t-s}$&$(z_k\psi_{l})_{t-s}$&$\delta_{k0}\varepsilon_{l,t-s}$&$(z_k\zeta_{l})_{t-s}$&$\delta_{k0}h_{l,t-s}$&$\delta_{k0}f_{l,t-s}$\\\hline
$\omega_k^{(s)}$&$0$&$0$&$\delta_{l0}\omega_{k,t-s}$&$0$&$\delta_{kl}\psi_{0,t-s}$&$0$&$0$&$0$\\\hline
$\theta_k^{(s)}$&$0$&$0$&$(z_l\theta_k)_{t-s}$&$(z_k\psi_l)_{t-s}$&$\delta_{k0}\beta(\varepsilon_{l,t-s})$&$0$&$0$&$\delta_{k0}\alpha(f_{l,t})$\\\hline
$f_k^{(s)}$&$\delta_{l0}\alpha(f_{k,t-s-1})$&$0$&$\delta_{l0}f_{k,t-s-1}$&$\begin{array}{ll}\delta_{l,h-3}(k+1)\cdot\\\theta_{h-3,t-s-1}\end{array}$&$0$&$\begin{array}{ll}\delta_{l,h-3}(k+1)\cdot\\z_{l,t-s}\end{array}$&$\delta_{kl}\psi_{0,t-s}$&$(M_\alpha)_{kl}\zeta_{0,t-s}$\\\hline
$h_k^{(s)}$&$0$&$0$&$\delta_{l0}h_{k,t-s-1}$&$0$&$0$&$\begin{array}{ll}\delta_{k,\frac{h-3}{2}}\delta_{l,h-3}\cdot\\\theta_{h-3,t-s}\end{array}$&$0$&$\delta_{kl}\psi_{0,t-s}$\\\hline
$\zeta_k^{(s)}$&$(z_l\psi_k)_{t-s-1}$&$0$&$(z_l\zeta_k)_{t-s-1}$&$\begin{array}{ll}\delta_{k,h-3}\delta_{l,h-3}\cdot\\\alpha(f_{\frac{h-3}{2},t-s-1})\end{array}$&$0$&$\begin{array}{ll}\delta_{k,h-3}\delta_{l,h-3}\cdot\\ f_{\frac{h-3}{2},t-s-1}\end{array}$&$\begin{array}{ll}\delta_{k,h-3}\delta_{l,\frac{h-3}{2}}\cdot\\\theta_{k,t-s}\end{array}$&$\begin{array}{ll}\delta_{k,h-3}(l+1)\cdot\\z_{k,t-s}\end{array}$\\\hline
$\varepsilon_k^{(s)}$&$-\delta_{l0}\beta(\varepsilon_{k,t-s})$&$\delta_{kl}\psi_{0,t-s-1}$&$\delta_{l,0}\varepsilon_{k,t-s-1}$&$0$&$\begin{array}{ll}-(M_\beta)_{k,l}\cdot\\\zeta_{0,t-s-1}\end{array}$&$0$&$0$&$0$\\\hline
$\psi_k^{(s)}$&$0$&$0$&$(z_k\psi_l)_{t-s-1}$&$0$&$0$&$\begin{array}{ll}\delta_{k,h-3}\delta_{l,h-3}\cdot\\\alpha(f_{h-3,t-s-1})\end{array}$&$0$&$\begin{array}{ll}\delta_{k,h-3}(l+1)\cdot\\\theta_{h-3,t-s}\end{array}$\\\hline
\end{tabular}
\caption{contraction map $\iota_a(b)$}
\label{contraction}
\end{table}

\end{landscape}

\begin{landscape}
\begin{table}
\begin{center}
\begin{tabular}[b]{|c||@{}c@{\hspace{0 cm}}|@{\hspace{-0.1 cm}}c@{\hspace{-0.18 cm}}|@{\hspace{-0.18 cm}}c@{\hspace{-0.18 cm}}|@{\hspace{-0.18 cm}}c@{\hspace{-0.18 cm}}|@{\hspace{-0.1 cm}}c@{\hspace{-0.18 cm}}|@{\hspace{-0.18 cm}}c@{\hspace{-0.18 cm}}|@{\hspace{-0.1 cm}}c@{\hspace{-0.18 cm}}|@{\hspace{-0.18 cm}}c@{\hspace{-0.18 cm}}|}
\hline
\backslashbox{$a$}{$b$}&$z_l^{(t)}$&$\omega_l^{(t)}$&$\theta_l^{(t)}$&$f_l^{(t)}$&$h_l^{(t)}$&$\zeta_l^{(t)}$&$\varepsilon_l^{(t)}$&$\psi_l^{(t)}$\\\hline\hline
$z_k^{(s)}$&$0$&$\begin{array}{cc}-\delta_{k0}sh\cdot\\\beta^{-1}(\omega_l^{(s+t)})\end{array}$&$\begin{array}{cc}(\frac{k}{2}-sh)\cdot\\(z_kz_l)^{(s+t)}\end{array}$&$0$&$\begin{array}{cc}-\delta_{k0}sh\cdot\\\alpha^{-1}(h_l^{(s+t)})\end{array}$&$0$&$0$&$\begin{array}{cc}(\frac{k}{2}-sh)\cdot\\(z_k\zeta_{l})^{(s+t)}\end{array}$\\\hline
$\omega_k^{(s)}$&&$0$&$0$&$0$&$0$&$0$&$                \begin{array}{cc}-(\frac{h}{2}+1+th)\cdot\\\delta_{kl}\zeta_0^{(s+t)}\end{array}  $&$0$\\\hline                                                                      $\theta_k^{(s)}$&&&$\begin{array}{cc}(\frac{l-k}{2}+(s-t)h)\cdot\\(z_k\theta_{l})^{(s+t)}\end{array}$&$\begin{array}{cc}-(1+th)\cdot\\\delta_{k0}f_l^{(s+t)}\end{array}$&$\begin{array}{cc}(-1+(s-t)h)\cdot\\\delta_{k0}h_l^{(s+t)}\end{array}$&$\begin{array}{cc}-(2+\frac{l}{2}+th)\cdot\\(z_k\zeta_{l})^{(s+t)}\end{array}$&$\begin{array}{cc}-(1+th+\frac{h}{2})\cdot\\\delta_{k0}\varepsilon_l^{(s+t)}\end{array}$&$\begin{array}{cc}-(2+\frac{k+l}{2}+(t-s)h)\cdot\\(z_k\psi_{l})^{(s+t)}\end{array}$\\\hline
$f_k^{(s)}$&&&&$0$&$\begin{array}{cc}-(1+sh)\cdot\\\delta_{kl}\zeta_0^{(s+t)}\end{array}$&$0$&$0$&$\begin{array}{cc}-(k+1)\cdot\\ (1+sh)\cdot\\\delta_{l,h-3}z_{h-3}^{(s+t+1)}\end{array}$\\\hline
$h_k^{(s)}$&&&&&$\begin{array}{c}(s-t)h\cdot\\(M_\alpha^{-1})_{kl}\cdot\\\psi_0^{(s+t)}\end{array}$&$\begin{array}{cc}-(\frac{h+1}{2}+th)\cdot\\\delta_{k,\frac{h-3}{2}}\delta_{l,h-3}\cdot\\z_{h-3}^{(s+t+1)}\end{array}$&$0$&$\begin{array}{c}((s-t)h-\frac{h-1}{2})\cdot\\\delta_{k,\frac{h-3}{2}}\delta_{l,h-3}\cdot\\\theta_{h-3}^{(s+t+1)}\end{array}$\\\hline
$\zeta_k^{(s)}$&&&&&&$0$&$0$&$\begin{array}{c}-(sh+\frac{h+1}{2})\cdot\\\delta_{k,h-3}\delta_{l,h-3}\cdot\\ f_{\frac{h-3}{2}}^{(s+t+1)}\end{array}$\\\hline
$\varepsilon_k^{(s)}$&&&&&&&$0$&$0$\\\hline
$\psi_k^{(s)}$&&&&&&&&$\begin{array}{c}(s-t)h\cdot\\\delta_{k,h-3}\delta_{l,h-3}\cdot\\\alpha(f_{\frac{h-3}{2}}^{(s+t+1)})\end{array}$\\\hline
\end{tabular}
\end{center}
\caption{Gerstenhaber bracket $[a,b]$}
\label{Gerstenhaber bracket}
\end{table}
\end{landscape}

\begin{landscape}
\thispagestyle{empty}	
 \begin{table}
 \begin{tabular}[b]{|@{\hspace{0.1 cm}}c@{\hspace{0.1 cm}}||@{\hspace{-0.18 cm}}c@{\hspace{-0.18 cm}}|@{\hspace{-0.18 cm}}c@{\hspace{-0.18 cm}}|@{\hspace{-0.18 cm}}c@{\hspace{-0.18 cm}}|@{\hspace{-0.18 cm}}c@{\hspace{-0.18 cm}}|@{\hspace{-0.18 cm}}c@{\hspace{-0.18 cm}}|@{\hspace{-0.18 cm}}c@{\hspace{-0.18 cm}}|@{\hspace{-0.18 cm}}c@{\hspace{-0.18 cm}}|@{\hspace{-0.18 cm}}c@{\hspace{-0.18 cm}}|}
\hline
\backslashbox{$a$}{$b$}&$\theta_{l,t}$&$\omega_{l,t}$&$z_{l,t}$&$\psi_{l,t}$&$\varepsilon_{l,t}$&$\zeta_{l,t}$&$h_{l,t}$&$f_{l,t}$\\\hline\hline
$\theta_k^{(s)}$&$\begin{array}{c}(1+\frac{l}{2}+th)\cdot\\(z_k\theta_l)_{t-s}\end{array}$&$\begin{array}{c}(\frac{1}{2}+t)h\cdot\\\delta_{k0}\omega_{l,t-s}\end{array}$&$\begin{array}{c}(1+\frac{k+l}{2}+(t-s)h)\cdot\\(z_kz_l)_{t-s}\end{array}$&$\begin{array}{c}((t+1)h-1-\frac{l}{2})\cdot\\(z_k\psi_l)_{t-s}\end{array}$&$\begin{array}{c}(\frac{1}{2}+(t-s))h\cdot\\\delta_{k0}\varepsilon_{l,t-s}\end{array}$&$\begin{array}{c}((t-s+1)h\\-1-\frac{l-k}{2})\\(z_k\zeta_l)_{t-s}\end{array}$&$\begin{array}{c}(t+1)h\cdot\\\delta_{k0}h_{l,t-s}\end{array}$&$\begin{array}{c}(t-s+1)h\cdot\\\delta_{k0}f_{l,t-s}\end{array}$\\\hline
$f_k^{(s)}$&$\begin{array}{c}-(1+sh)\cdot\\\delta_{l0}f_{k,t-s-1}\end{array}$&$0$&$0$&$\begin{array}{c}-(1+sh)\cdot\\\delta_{l,h-3}z_{h-3,t-s}\end{array}$&$0$&$0$&$\begin{array}{c}-(1+sh)\cdot\\\delta_{kl}\zeta_{0,t-s}\end{array}$&$0$\\\hline
$h_k^{(s)}$&$\begin{array}{c}(1+th)\cdot\\\delta_{l0}h_{k,t-s-1}\end{array}$&$0$&$\begin{array}{c}\delta_{l0}(t-s)h\cdot\\\alpha^{-1}(h_{k,t-s-1})\end{array}$&$\begin{array}{c}(th+\frac{h+1}{2})\cdot\\\delta_{k,\frac{h-3}{2}}\delta_{l,h-3}\cdot\\\theta_{h-3,t-s}\end{array}$&$0$&$\begin{array}{c}((t-s)h+\frac{h-1}{2})\cdot\\\delta_{k,\frac{h-3}{2}}\delta_{l,h-3}\cdot\\z_{h-3,t-s}\end{array}$&$\begin{array}{c}(t+1)h\cdot\\(M_\alpha^{-1})_{lk}\cdot\\\psi_{0,t-s}\end{array}$&$\begin{array}{c}((t-s+1)h-1)\cdot\\\delta_{kl}\zeta_{0,t-s}\end{array}$\\\hline
$\zeta_k^{(s)}$&$\begin{array}{c}-(2+\frac{k}{2}+sh)\cdot\\(z_l\zeta_k)_{t-s-1}\end{array}$&$0$&$0$&$\begin{array}{c}-(sh+\frac{h+1}{2})\cdot\\\delta_{k,h-3}\delta_{l,h-3}\cdot\\f_{\frac{h-3}{2},t-s-1}\end{array}$&$0$&$0$&$\begin{array}{c}-(sh+\frac{h+1}{2})\cdot\\\delta_{k,h-3}\delta_{l,\frac{h-3}{2}}\cdot\\z_{h-3,t-s}\end{array}$&$0$\\\hline
$\varepsilon_k^{(s)}$&$\begin{array}{c}((s+\frac{1}{2})h+1)\cdot\\\delta_{l0}\varepsilon_{k,t-s-1}\end{array}$&$\begin{array}{c}-((s+\frac{1}{2})h+1)\cdot\\\delta_{kl}\zeta_{0,t-s-1}\end{array}$&$0$&$0$&$0$&$0$&$0$&$0$\\\hline
$\psi_k^{(s)}$&$\begin{array}{c}(1+\frac{l}{2}+th)\cdot\\(z_l\psi_k)_{t-s-1}\end{array}$&$0$&$\begin{array}{c}((t-s)h\\-1-\frac{k-l}{2})\cdot\\(z_l\zeta_k)_{t-s-1}\end{array}$&$\begin{array}{c}(th+\frac{h+1}{2})\cdot\\\delta_{k,h-3}\delta_{l,h-3}\cdot\\\alpha(f_{\frac{h-3}{2},t-s-1})\end{array}$&$0$&$\begin{array}{c}(t-s)h\cdot\\\delta_{k,h-3}\delta_{l,h-3}\\f_{\frac{h-3}{2},t-s-1}\end{array}$&$\begin{array}{c}(t+1)h\cdot\\\delta_{k,h-3}\delta_{l,h-3}\cdot\\\theta_{h-3,t-s}\end{array}$&$\begin{array}{c}(l+1)\cdot\\
((t-s)h\\+1+\frac{h-3}{2})\cdot\\\delta_{k,h-3}\\z_{h-3,t-s}\end{array}$\\\hline
$z_k^{(s)}$&$\begin{array}{c}(k-sh)\\(z_k\theta_l)_{t-s}\end{array}$&$\begin{array}{c}-\delta_{k0}sh\\\beta^{-1}(\omega_{l,t-s})\end{array}$&$0$&$\begin{array}{c}(\frac{k}{2}-sh)\cdot\\(z_k\zeta_l)_{t-s}\end{array}$&$0$&$0$&$\begin{array}{c}(k-sh)\cdot\\\alpha^{-1}(h_{l,t-s})\end{array}$&$0$\\\hline
$\omega_{k}^{(s)}$&$\begin{array}{c}(1+th)\cdot\\\delta_{l0}\omega_{k,t-s}\end{array}$&$0$&$\begin{array}{c}\delta_{l0}(\frac{1}{2}+t-s)h\cdot\\\beta^{-1}(\omega_{k,t-s})\end{array}$&$0$&$\begin{array}{c}\delta_{kl}\cdot(-1+h\\+(t-s)h)\cdot\\\zeta_{0,t-s}\end{array}$&$0$&$0$&$0$\\\hline
\end{tabular}
\caption{Lie derivative $\mathcal{L}_a(b)$}
\label{Lie derivative}
\end{table}
\end{landscape}

\section{Batalin-Vilkovisky structure on Hochschild cohomology}
Recall the isomorphism (\ref{H-dual}). It translates the Connes differential $B:HH_\bullet(A)\rightarrow HH_{\bullet+1}(A)$ on Hochschild homology into a differential $\Delta:HH^\bullet(A)\rightarrow HH^{\bullet -1}(A)$ on Hochschild cohomology, i.e. we have the commutative diagram

$$
\CD
HH_\bullet(A)     @>B>>  HH_{\bullet+1}(A)\\
@V\mathbb{D}V\sim V      @V\sim V\mathbb{D}V\\         
HH^{6m+2-\bullet}(A)[2mh+2] @>\Delta>>  HH_{6m+1-\bullet}(A)[2mh+2]
\endCD
$$

\begin{theorem}\emph{(BV structure on Hochschild cohomology)}
 $\Delta$ makes $HH^\bullet(A)$ a Batalin-Vilkovisky algebra, i.e. 
 for the Gerstenhaber bracket we get the following equation:
\begin{equation}\label{BV-identity}
 [a,b]=\Delta(a\cup b)-\Delta(a)\cup b-(-1)^{|a|}a\cup\Delta(b),\qquad\forall a,b\in HH^*(A).
\end{equation}

The isomorphism $\mathbb{D}$ intertwines contraction and cup-product maps, i.e. we have
\begin{equation}\label{intertwining}
 \mathbb{D}(\iota_\eta c)=\eta\cup\mathbb{D}(c),\qquad\forall c\in HH_\bullet(A),\,\eta\in HH^\bullet(A).
\end{equation}
\end{theorem}
\begin{remark}
 Note that $\Delta$ in equation (\ref{BV-identity}) depends on which $m\in\mathbb{N}$ we choose to identify $\mathbb{D}:HH_\bullet(A)\stackrel{\sim}{\rightarrow} HH^{6m+2-\bullet}(A)[2mh+2]$, where the Gerstenhaber bracket does not.\\
\end{remark}

\begin{proof}
We apply the functor 
\[
\begin{array}{rcl}
 \mathrm{Hom}_{A^e}(-,A\otimes_\mathbb{C} A):A^e-\mathrm{mod}&\rightarrow&A^e-\mathrm{mod},\\
 M&\mapsto&M^\vee
 \end{array}
\]
on the Schofield resolution:
\begin{equation}\label{Hom-S}
 \begin{array}{rcl}
(A\otimes A)^\vee&\stackrel{d_1^\vee}{\rightarrow}&(A\otimes V\otimes A)^\vee\stackrel{d_2^\vee}{\rightarrow}(A\otimes A[2])^\vee\stackrel{d_3^\vee}{\rightarrow}(A\otimes\mathfrak{N}[h])^\vee\stackrel{d_4^\vee}{\rightarrow}\\
&\stackrel{d_4^\vee}{\rightarrow}&(A\otimes V\otimes\mathfrak{N}[h])^\vee\stackrel{d_5^\vee}{\rightarrow}(A\otimes\mathfrak{N}[h+2])^\vee\stackrel{d_6^\vee}{\rightarrow}(A\otimes A[2h])^\vee\stackrel{d_7^\vee}{\rightarrow}\ldots
 \end{array}
 \end{equation}
An element in $(A\otimes A)^\vee$ or $(A\otimes\mathfrak{N})^\vee$ is determined by the image of $1\otimes1$,\\
An element in $(A\otimes V\otimes A)^\vee$ or $(A\otimes V\otimes\mathfrak{N})^\vee$ by the images of $1\otimes a\otimes1$ for all arrows $a\in\bar Q$.

Let us define $\sigma=\left\{\begin{array}{cc}+1&Q=A,\\-1&Q=D,E\end{array}\right.$.

We make the following identifications:

\underline{$(A\otimes A)[-2mh]=(A\otimes A[2mh])^\vee$:}\\
we identify $x\otimes y$ with the map that sends $1\otimes 1$ to $\sigma^m y\otimes x$,

\underline{$(A\otimes V\otimes A)[-2mh-2]=(A\otimes V\otimes A[2mh])^\vee$:}\\
we identify $\sum\limits_{a\in\bar Q}\epsilon_ax_a\otimes a^*\otimes y_a$ with the map that sends $1\otimes a\otimes 1$ to $-\sigma^my_a\otimes x_a$,

\underline{$(A\otimes A)[-2mh-2]=(A\otimes A[2mh+2])^\vee$}:\\
we identify $x\otimes y$ with the map that sends $1\otimes1$ to $-\sigma^my\otimes x$,

\underline{$(A\otimes \mathfrak{N})[-(2m+1)h]=(A\otimes \mathfrak{N}[(2m+1)h])^\vee$:} \\
we identify $x\otimes y$ with the map that sends $1\otimes 1$ to $-\sigma^m\eta(y)\otimes x$,

\underline{$(A\otimes V\otimes \mathfrak{N})[-(2m+1)h-2]=(A\otimes V\otimes \mathfrak{N}[(2m+1)h])^\vee$:}\\
we identify $\sum\limits_{a\in\bar Q}\epsilon_ax_a\otimes \eta(a^*)\otimes y_a$ with the map that sends $1\otimes a\otimes 1$ to $\sigma^{m+1} \eta(y_a)\otimes x_a$,

\underline{$(A\otimes \mathfrak{N})[-2(m+1)h-2]=(A\otimes A[2(m+1)h+2])^\vee$}:\\
we identify $x\otimes y$ with the map that sends $1\otimes1$ to $\sigma^{m+1}\eta(y)\otimes x$,

so (\ref{Hom-S}) becomes

\begin{equation} \label{Hom-S2}
 \begin{array}{rcl}
(A\otimes A)&\stackrel{d_1^\vee}{\rightarrow}&(A\otimes V\otimes A)\stackrel{d_2^\vee}{\rightarrow}(A\otimes A[-2])\stackrel{d_3^\vee}{\rightarrow}(A\otimes\mathfrak{N}[-h])\stackrel{d_4^\vee}{\rightarrow}\\
&\stackrel{d_4^\vee}{\rightarrow}&(A\otimes V\otimes\mathfrak{N}[-h-2])\stackrel{d_5^\vee}{\rightarrow}(A\otimes\mathfrak{N}[-h-2])\stackrel{d_6^\vee}{\rightarrow}(A\otimes A[-2h])\stackrel{d_7^\vee}{\rightarrow}\ldots
 \end{array}
\end{equation}
We show under the identification from above, the differentials $d_i^\vee$ corresponds to the differentials from the Schofield resolution, i.e. (\ref{Hom-S2})
can be rewritten in this form:

\begin{equation}
 \begin{array}{rcl}
(A\otimes A)&\stackrel{d_2[-2]}{\rightarrow}&(A\otimes V\otimes A)\stackrel{d_1[-2]}{\rightarrow}(A\otimes A[-2])\stackrel{d_6[-2h-2]}{\rightarrow}\\
&\stackrel{d_6[-2h-2]}{\rightarrow}&(A\otimes\mathfrak{N}[-h])\stackrel{d_5[-2h-2]}{\rightarrow}(A\otimes V\otimes\mathfrak{N}[-h-2])\stackrel{d_4[-2h-2]}{\rightarrow}\\
&\stackrel{d_4[-2h-2]}{\rightarrow}&(A\otimes\mathfrak{N}[-h-2])\stackrel{d_3[-2h-2]}{\rightarrow}(A\otimes A[-2h])\stackrel{d_2[-2h-2]}{\rightarrow}\ldots
 \end{array}
\end{equation}
It is enough to show this for the first period.
\[
 d_1^\vee(x\otimes y)(1\otimes a\otimes1)=(x\otimes y)\circ(a\otimes1-1\otimes a)=ay\otimes x-y\otimes xa,
\]
so
\[
 d_1^\vee(x\otimes y)=\sum\limits_{a\in\bar Q}\epsilon_a(xa\otimes a^*\otimes y-x\otimes a^*\otimes ay)=\sum\limits_{a\in\bar Q}\epsilon_a(xa\otimes a^*\otimes y+x\otimes a\otimes a^*y)=d_2(x\otimes y),
\]

\begin{eqnarray*}
 d_2^\vee(\sum\limits_{a\in\bar Q}\epsilon_ax_a\otimes a^*\otimes y_a)(1\otimes1)&=&(\sum\limits_{a\in\bar Q}\epsilon_ax_a\otimes a^*\otimes y_a)\circ(\sum\limits_{b\in\bar Q}\epsilon_bb\otimes b^*\otimes1+\epsilon_b1\otimes b\otimes b^*)\\
 &=&\sum\limits_{a\in\bar Q}(\epsilon_aa^*y_a\otimes x_a-\epsilon_ay_a\otimes x_aa^*),
\end{eqnarray*}
so
\[
 d_2^\vee(\sum\limits_{a\in\bar Q}\epsilon_ax_a\otimes a^*\otimes y_a)=\sum\limits_{a\in\bar Q}\epsilon_a(x_aa^*\otimes y_a-x_a\otimes a^*y_a)=d_1(\epsilon_ax_a\otimes a^*\otimes y_a),
\]
\[
 d_3^\vee(x\otimes y)(1\otimes1)=(x\otimes y)\circ(\sum\limits_{x_i\in B}x_i\otimes x_i^*)=-\sum\limits_{x_i\in B}x_iy\otimes xx_i^*=-\sum\limits_{x_i\in B}\eta(x_i^*)y\otimes xx_i,
\]
so
\[
d_3^\vee(x\otimes y)=\sum\limits_{x_i\in B}xx_i\otimes x_i^*\eta(y)=\sum\limits_{x_i\in B}xyx_i\otimes x_i^*=d_6(x\otimes y)
\]

\[
 d_4^\vee(x\otimes y)(1\otimes a\otimes 1)=(x\otimes y)\circ(a\otimes 1-1\otimes a)=-a\eta(y)\otimes x+\eta(y)\otimes x\eta(a),
\]
so
\begin{eqnarray*}
 d_4^\vee(x\otimes y)&=&\sum\limits_{a\in\bar Q}\epsilon_a\sigma(-x\otimes\eta(a^*)\otimes\eta(a)y+x\eta(a)\otimes\eta(a^*)\otimes y)\\
 &=&\sum\limits_{a\in\bar Q}(\epsilon_axa\otimes a^*\otimes y+\epsilon_ax\otimes a\otimes a^*y)=d_5(x\otimes y),
\end{eqnarray*}

\[
\begin{array}{l}
 d_5^\vee(\sum\limits_{a\in\bar Q}\epsilon_ax_a\otimes \eta(a^*)\otimes y_a)(1\otimes a\otimes1)\\
 \qquad=(\sum\limits_{a\in\bar Q}\epsilon_ax_a\otimes \eta(a^*)\otimes y_a)\circ(\sum\limits_{b\in\bar Q}(\epsilon_bb\otimes b^*\otimes1+\epsilon_b1\otimes b\otimes b^*))\\
 \qquad=\sigma\sum\limits_{a\in\bar Q}(-\epsilon_a^*\eta(y_a)\otimes x_a+\epsilon_a\eta(y_a)\otimes x_a\eta(a^*)),
\end{array}
\]
so

\begin{eqnarray*}
d_5^\vee(\sum\limits_{a\in\bar Q}\epsilon_ax_a\otimes \eta(a^*)\otimes y_a&=&\sum\limits_{a\in\bar Q}(-\epsilon_ax_a\otimes\eta(a^*)y_a+\epsilon_ax_a\eta(a^*)\otimes y_a)\\
&=&d_4(\sum\limits_{a\in\bar Q}\epsilon_ax_a\otimes\eta(a^*)\otimes y_a),
\end{eqnarray*}

\begin{eqnarray*}
 d_6^\vee(x\otimes y)(1\otimes1)&=&(x\otimes y)\circ(\sum\limits_{x_i\in B}x_\otimes x_i^*)=\sigma\sum\limits_{x_i\in B}x_i\eta(y)\otimes x\eta(x_i^*)\\&=&\sigma\sum\limits_{x_i\in B}x_i^*\eta(y)\otimes xx_i=\sigma\sum\limits_{x_i\in B}x_i^*\otimes xyx_i,
\end{eqnarray*}
so
\[
 d_6^\vee(x\otimes y)=\sum\limits_{x_i\in B}xyx_i\otimes x_i^*=d_3(x\otimes y)
\]

Fix $m\geq0$. The map which shifts the degree by $-2mh-2$ produces the following diagram which commutes by the computations above:

$$
\CD
A\otimes A[2mh+2]@>d_{6m+2}>>A\otimes V\otimes A[2h]@>d_{6m+1}>>\ldots\\
@VVV @VVV \\
(A\otimes A)@>d_1^\vee>>(A\otimes V\otimes A)[-2]@>d_2^\vee>>\ldots
\endCD
$$

$$
\CD
\ldots@>d_2>>A\otimes V\otimes A@>d_1>>A\otimes A@>mult.>>A\\
@.@VVV@VVV\\
\ldots@>d_{6m+1}^\vee>>(A\otimes V\otimes A)[-2mh-2]@>d_{6m+2}^\vee>>(A\otimes A)[-mh-2]@>mult.>>A
\endCD
$$

Similarly to the proof of \cite[Theorem 3.4.3.]{G}, this self-dual morphism of the Schofield resolution $C^\bullet$ into the dual complex $(C^\bullet)^\vee$ can be used to prove (\ref{intertwining}).

(\ref{BV-identity}) follows, as in the proof of \cite[Theorem 3.4.3.]{G}, from (\ref{intertwining}) and the calculus structure. 

\end{proof}

\subsection{Computation of the calculus structure of the preprojective algebra}
Since the calculus structure is defined on Hochschild chains and cochains, we have to work with the on the resolution for computations. It turns out that we only have to compute $\mathfrak{L}_{\theta_0}$ directly, the rest can be deduced from formulas given by the calculus and the BV structure.

$$
\CD
\hdots@>d_3>>A\otimes A[2]@>d_2>>A\otimes V\otimes A@>d_1>>A\otimes A@>d_0>> A @>>> 0\\
@. @V\mu_2 VV @V\mu_1 VV @\vert @\vert\\
\hdots@>b_3>>A^{\otimes4}@>b_2>>A^{\otimes3}@>b_1>>A^{\otimes2}@>b_0>>A@>>>0
\endCD
$$
These maps $\psi_i$ gives us a chain map between the Schofield and the bar resolution:
\begin{eqnarray*}
\mu_1(1\otimes y\otimes 1)&=&1\otimes y\otimes 1,\\
\mu_2(1\otimes 1)&=&\sum\limits_{a\in\bar Q}\epsilon_a 1\otimes a\otimes a^*\otimes1,\\
\mu_3(1\otimes1)&=&\sum\limits_{a\in\bar Q}\sum\limits_{x_i\in B}\epsilon_a 1\otimes x_i\otimes a\otimes a^*\otimes x_i^*,
\end{eqnarray*}
and 
\[
 \mu_{3+i}=\mu_i\sum\limits_{a\in\bar Q}\sum\limits_{x_i\in B}x_i\otimes a\otimes a^*\otimes x_i^*.
\]

Now, we apply the functor $-\otimes_{A^e} A$ on the commutative diagram:

$$
\CD
\hdots@>d_3'>>A^R[2]@>d_2'>>(V\otimes A)^R @>d_1>>A^R@>>> 0 \\
@. @V\mu_2' VV @V\mu_1' VV @\vert \\
\hdots@>b_3'>>(A^{\otimes3})^R@>b_2>>(A^{\otimes2})^R@>b_1>>(A^{\otimes1})^R@>>>0
\endCD
$$

where

\begin{eqnarray*}
\mu_1'(x\otimes y)&=&x\otimes y,\\
\mu_2'(x)&=&\sum\limits_{a\in\bar Q}\epsilon_a a\otimes a^*\otimes x,\\
\mu_3'(x)&=&\sum\limits_{a\in\bar Q}\sum\limits_{x_i\in B}\epsilon_a x_i\otimes a\otimes a^*\otimes x_i^*\eta(x),
\end{eqnarray*}
and

\[
 \mu_{3+i}'=\mu_i'\sum\limits_{a\in\bar Q}\sum\limits_{x_i\in B}x_i\otimes a\otimes a^*\otimes x_i^*.
\]

Now, we compute $\mathfrak{L}_{\theta_0}$:

\begin{lemma}\label{Lie theta}
 For each $x\in HH_i(A)$, 
 \begin{equation}
  \mathfrak{L}_{\theta_0}(x)=x\frac{\deg(x)}{2}.
 \end{equation}
\end{lemma}

\begin{proof}
 Via $\mu'$, we already identified $x\in HH_i(A)$ with cycles in the Hochschild chain, but we still have to identify $\theta_0$ with an element in $\mathrm{Hom}_{A^e}(A^{\otimes3},A)$:

given any monomial $b=b_1\ldots b_l$, $b_i\in V$, the map 
\[
\tau(1\otimes b\otimes1)=\sum\limits_{i=1}^lb_1\ldots b_{i-1}\otimes b_i\otimes b_{i+1}\ldots b_l
\] 
makes the diagram
$$
\CD
A\otimes V\otimes A@>d_1>>A\otimes A@>d_0>> A @>>> 0\\
@A\tau AA @\vert @\vert\\
A^{\otimes3}@>b_1>>A^{\otimes2}@>b_0>>A@>>>0
\endCD
$$
commute.

Applying $\mathrm{Hom}_{A^e}(\_\_\otimes A)$, we get a map

\[
 \tau^*: \mathrm{Hom}_k(V)\rightarrow\mathrm{Hom}_k(A),
\]
such that
\[
 (\theta_0\circ\tau^*)(b_1\ldots b_l)=\sum\limits_{i=1}^lb_1\ldots b_{i-1}\theta_0(b_i) b_{i+1}\ldots b_l=s(b)\cdot b,
\]
where for $b=b_1\ldots b_l$, $s(b)$ is the number of $b_i\in Q^*$.

Recall from \cite[(3.5), page 46]{D} that the Lie derivative of $\theta_0\circ\tau^*$ on Hochschild chains is defined by

\begin{eqnarray*}
 \mathfrak{L}_{\theta_0\circ\tau^*}(a_1\otimes\cdots\otimes a_k)&=&\sum\limits_{i=1}^k a_1\otimes\cdots\otimes (\theta_0\circ\tau^*)(a_i)\otimes\cdots\otimes a_k\\
 &=&\sum\limits_{i=1}^k (s(a_1)+\cdots s(a_k)) a_1\otimes\cdots\otimes a_k,
\end{eqnarray*}
and it can easily be checked that for each $x\in HH_i(A)$, $\mathfrak{L}_{\theta_0\circ\tau^*}$ acts on $\mu_i'(x)$, $x\in HH^i(A)$, by multiplication with $\frac{1}{2}\deg(x)$.

\end{proof}

\newpage

\subsubsection{The contraction map}
From (\ref{intertwining}) we know that the contraction map on Hochschild homology is given by the cup product on Hochschild cohomology which was computed in \cite{ES2} and \cite{Eu2}. Table \ref{contraction} contains these results, rewritten in terms of the contraction maps.

\subsubsection{The Connes differential}
We start with the computation of the Connes differential and recall the diagram from \cite{EE2}:

$$
\CD
\text{degree}\\
@.                                0                    \\
@.                              @VVV\\
 2\leq\deg\leq h-1   @.  {HH_1}(A)@=   U     \\
@.                              @VB_1VV    @V\sim VV\\
 2\leq\deg\leq h     @.  {HH_2}(A)@=   U     @. \oplus @.    Y[h] \\
@.                              @VB_2VV               @.    @.    @V\sim VV\\
 h\leq\deg\leq 2h-2  @.  {HH_3}(A)@= U^*[2h] @. \oplus @. Y^*[h] \\
@.                              @VB_3VV     @V\sim VV\\
h+1\leq\deg\leq2h-2  @.  {HH_4}(A)@= U^*[2h] \\
@.                              @VB_4VV       @V0VV \\
2h    @.  {HH_5}(A)@=   K^*[2h]  \\
@.                              @VB_5VV     @V\sim VV\\
2h                   @.  {HH_6}(A)@=   K[2h]  \\
@.                              @VB_6VV       @V0VV    \\
2h+2\leq\deg\leq3h-1\, @.  {HH_7}(A)@=U[2h]\\
@.                              @VB_7VV\\
@.                              \vdots
\endCD
$$

\begin{proposition}
 The Connes differential $B$ is given as follows:
 \begin{eqnarray*}
 B_{1+6s}(\theta_k^{(s)})&=&(1+\frac{k}{2}+sh)z_k^{(s)},\\
 B_{2+6s}(\omega_k^{(s)})&=&(\frac{1}{2}+s)h\beta^{-1}(\omega_k^{(s)}),\\
 B_{2+6s}(z_{k}^{(s)})&=&0,\\
 B_{3+6s}(\psi_k^{(s)})&=&((s+1)h-1-\frac{k}{2})\zeta_k,\\
 B_{3+6s}(\varepsilon_k^{(s)})&=&0,\\
 B_{4+6s}&=&0,\\
 B_{5+6s}(h_k^{(s)})&=&(s+1)h\alpha^{-1}(h_k^{(s)}),\\
 B_{6+6s}&=&0.
 \end{eqnarray*}
\end{proposition}
\begin{proof}
We use the Cartan identity (\ref{Cartan}) with $a\in\theta_0$,
\begin{equation}
\mathcal{L}_{\theta_0}=B\iota_{\theta_0}+\iota_{\theta_0} B,
\end{equation}
where $\mathcal{L}_{\theta_0}$ acts on $x\in HH_i$ by multiplication by $\frac{1}{2}\deg(x)$ (see Lemma (\ref{Lie theta})).
The above identities for the Connes differential follow since $\iota_{\theta_0}$ acts on $\theta_k^{(t)}$, $\omega_k^{(t)}$, $\psi_k^{(t)}$ and $h_k^{(t)}$ by zero, and $z_k^{(t)}$, $\beta^{-1}(\omega_k^{(t)})$, $\zeta_k^{(t)}$ and $\alpha^{-1}(h_k^{(s)})$ are their unique preimages the contraction with $\iota_{\theta_0}$.
\end{proof}

\subsubsection{The Gerstenhaber bracket}

We compute the brackets using the identification \\$HH^{i}(A)=HH_{6m+2-i}(A)[-2mh-2]$ for $m>>1$ and the BV-identity (\ref{BV-identity}).

\textbf{\underline{Brackets involving $HH^{6s}(A)$}}:\\

By degree argument these brackets are zero:\\
$\omega_k^{(s)}$ with $HH^{1+6t}(A)$,  $HH^{2+6t}(A)$, $HH^{3+6t}(A)$, $HH^{4+6t}(A)$, and $\psi_l\in HH^{5+6t}(A)$. 
From the BV-identity (\ref{BV-identity}), we see that brackets of $z_k^{(s)}$ with $z_l^{(t)}\in HH^{6t}(A)$, $HH^{2+6t}(A)$, $HH^{4+6t}(A)$ and $\varepsilon_l^{(t)}\in HH^{5+6t}(A)$ are zero because $\Delta$ acts by zero on $U[-2th-2]\subset HH^{6t}(A)$, $HH^{2+6t}(A)$ and $HH^{4+6t}(A)$.

We compute the remaining brackets:
\begin{eqnarray*}
 [z_k^{(s)},\omega_l^{(t)}]&=&\Delta(z_k^{(s)}\cup\omega_l^{(t)})-\underbrace{\Delta(z_k^{(s)})}_{=0}\cup\omega_l^{(t)}-z_k^{(s)}\cup\Delta(\omega_l^{(t)})\\
 &=&\delta_{k0}\Delta(\omega_l^{(s+t)})-(\frac{h}{2}+(m-t)h)z_k^{(s)}\cup\beta^{-1}(\omega_l^{(t)})\\
 &=&\delta_{k0}(\frac{h}{2}+(m-s-t)h)\beta^{-1}(\omega_l^{(s+t)})-(\frac{h}{2}+(m-t)h)\beta^{-1}(\omega_l^{(s+t)})\\
 &=&-\delta_{k0}sh\beta^{-1}(\omega_l^{(s+t)}),\\
 {}[z_k^{(s)},\theta_l^{(t)}]&=&\Delta(z_k^{(s)}\cup\theta_l^{(t)})-\underbrace{\Delta(z_k^{(s)})}_{=0}\cup\theta_l^{(t)}-z_k^{(s)}\cup\Delta(\theta_l^{(t)})\\
 &=&\Delta((z_k\theta_{l})^{(s+t)})-(1+\frac{l}{2}+(m-t)h)z_k^{(s)}z_l^{(t)}\\
 &=&(1+\frac{k+l}{2}+(m-s-t)h)(z_{k}z_l)^{(s+t)}\\&&-(1+\frac{l}{2}+(m-t)h)(z_{k}z_l)^{(s+t)}\\
 &=&(\frac{k}{2}-sh)(z_kz_l)^{(s+t)},\\
\end{eqnarray*}
\begin{eqnarray*}
[z_k^{(s)},h_l^{(t)}] &=&\Delta(z_k^{(s)}\cup h_l^{(t)})-\underbrace{\Delta(z_k^{(s)})}_{=0}\cup h_l^{(t)}-z_k^{(s)}\cup\Delta(h_l^{(t)})\\
&=&\delta_{k0}\Delta(h_l^{(s+t)})-(h+(m-t-1)h)z_k^{(s)}\cup\alpha^{-1} (h_l^{(t)})\\
&=&\delta_{k0}(h+(m-s-t-1)h)\alpha^{-1}(h_l^{(s+t)})\\&&-\delta_{k0}(h+(m-t-1)h)\alpha^{-1}(h_l^{(s+t)})\\
&=&-\delta_{k0}sh\alpha^{-1}(h_l^{(s+t)}),
\end{eqnarray*}
\begin{eqnarray*}
 [z_k^{(s)},\psi_l^{(t)}]&=&\Delta(z_k^{(s)}\cup\psi_l^{(t)})-\underbrace{\Delta(z_k^{(s)})}_{=0}\cup\psi_l^{(t)}-z_k^{(s)}\cup\Delta(\psi_l^{(t)})\\
 &=&\Delta((z_k\psi_{l})^{(s+t)})-(h-1-\frac{t}{2})z_k^{(s)}\zeta_l^{(t)}\\
 &=&((m-s-t)h-1-\frac{l-k}{2})(z_k\zeta_{l})^{(s+t)}\\&&-((m-t)h-1-\frac{l}{2})(z_k\zeta_{l})^{(s+t)}\\
 &=&(\frac{k}{2}-sh)(z_k\zeta_{l})^{(s+t)}
\end{eqnarray*}

\begin{eqnarray*}
 [\omega_k^{(s)},\varepsilon_l^{(t)}]&=&\Delta(\omega_k^{(s)}\cup\varepsilon_l^{(t)})-\Delta(\omega_k^{(s)})\cup\varepsilon_l^{(t)}-\omega_k^{(s)}\cup\underbrace{\Delta\varepsilon_l^{(t)}}_{=0}\\
 &=&\Delta(\delta_{kl}\psi_0^{{s+t}})-(\frac{h}{2}+(m-s)h)\beta^{-1}(\omega_k^{(s)})\cup\varepsilon_l^{(t)}\\
 &=&\delta_{kl}(h-1+(m-s-t-1)h)\zeta_0^{(s+t)}-\delta_{kl}(\frac{h}{2}+(m-s)h)\zeta_0\\
 &=&\delta_{kl}(-\frac{h}{2}-1-th)\zeta_0
\end{eqnarray*}

\textbf{\underline{Brackets involving $HH^{1+6s}(A)$}}:\\

\begin{eqnarray*}
  [\theta_k^{(s)},\theta_l^{(t)}]&=&\Delta(\underbrace{\theta_k^{(s)}\cup\theta_l^{(t)}}_{=0})-\Delta(\theta_k^{(s)})\cup\theta_l^{(t)}+\theta_k^{(s)}\cup\Delta(\theta_l^{(t)})\\
  &=&-(1+\frac{k}{2}+(m-s)h)z_k^{(s)}\theta_l^{(t)}+(1+\frac{l}{2}+(m-t)h))\theta_k^{(s)}z_l^{(t)}\\
  &=&(\frac{l-k}{2}+(s-t)h)(z_k\theta_{l})^{(s+t)}
\end{eqnarray*}
\begin{eqnarray*}
[\theta_k^{(s)},f_l^{(t)}]&=&\Delta(\theta_k^{(s)}\cup f_l^{(t)})-\Delta(\theta_k^{(s)})\cup f_l^{(t)}+\theta_k^{(s)}\cup\underbrace{\Delta(f_l^{(t)})}_{=0}\\
&=&\delta_{k0}(\Delta(\alpha(f_l^{(s+t)}))-(1+(m-s)h)f_l^{(s+t)})\\
&=&\delta_{k0}(h+(m-s-t-1)h)f_l^{(s+t)}-(1+(m-s)h)f_l^{(s+t)}\\
&=&-\delta_{k0}(1+th)f_l^{(s+t)}
\end{eqnarray*}
\begin{eqnarray*}
[\theta_k^{(s)},h_l^{(t)}]&=&\Delta(\underbrace{\theta_k^{(s)}\cup h_l^{(t)}}_{=0})-\Delta(\theta_k^{(s)})\cup h_l^{(t)}+\theta_k^{(s)}\cup\Delta(h_l^{(t)})\\
&=&-(1+(m-s)h+\frac{k}{2})z_k^{(s)}\cup h_l^{(t)}+(h+(m-t-1)h)\theta_k^{(s)}\cup\alpha^{-1}(h_l^{(t)})\\
&=&-\delta_{k0}(1+(m-s)h)h_l^{(s+t)}\\&&+\delta_{k0}(m-t)hh_l^{(s+t)}\\
&=&\delta_{k0}(-1+(s-t)h)h_l^{(s+t)}
\end{eqnarray*}
\begin{eqnarray*}
 [\theta_k^{(s)},\zeta_l^{(t)}]&=&\Delta(\theta_k^{(s)}\cup\zeta_l^{(t)})-\Delta(\theta_k^{(s)})\cup\zeta_l^{(t)}+\theta_k^{(s)}\cup\underbrace{\Delta(\zeta_l^{(t)})}_{=0}\\
 &=&\Delta((z_k\psi_{l})^{(s+t)})-(1+\frac{k}{2}+(m-s)h)z_k^{(s)}\cup\zeta_l^{(t)}\\
 &=&(h-1-\frac{l-k}{2}+(m-s-t-1)h)(z_k\zeta_{l})^{(s+t)}\\&&-(1+\frac{k}{2}+(m-s)h)(z_k\zeta_{l})^{(s+t)}\\
 &=&-(2+\frac{l}{2}+th)(z_k\zeta_{l})^{(s+t)}
\end{eqnarray*}
\begin{eqnarray*}
[\theta_k^{(s)},\psi_l^{(t)}]&=&\underbrace{\Delta(\theta_k^{(s)}\cup\psi_l^{(t)})}_{=0}-\Delta(\theta_k^{(s)})\cup\psi_{l}^{(t)}+\theta_{k}^{(s)}\cup\Delta(\psi_l^{(t)})\\
&=&-(1+\frac{k}{2}+(m-s)h)z_k^{(s)}\psi_{l}^{(t)}+(h-1-\frac{l}{2}+(m-t-1)h)\theta_{k}^{(s)}\zeta_{l}^{(t)}\\
&=&-(2+\frac{k+l}{2}+(t-s)h)(z_k\psi_{l})^{(s+t)},
\end{eqnarray*}

\begin{eqnarray*}
[\theta_k^{(s)},\varepsilon_{l}^{(t)}]&=&\Delta(\theta_k^{(s)}\cup\varepsilon_{l}^{(t)})-\Delta(\theta_k^{(s)})\cup\varepsilon_l^{(t)}+\theta_k^{(s)}\cup\underbrace{\Delta(\varepsilon_l^{(t)})}_{=0}\\
&=&\delta_{k0}\Delta(\beta(\varepsilon_l^{(s+t)}))-\delta_{k0}(1+(m-s)h+\frac{k}{2})z_k^{(s)}\varepsilon_l^{(t)}\\
&=&\delta_{k0}(\frac{h}{2}+(m-s-t-1)h)\varepsilon_l^{(s+t)}-(1+(m-s)h)\varepsilon_l^{(s+t)}\\
&=&-\delta_{k0}(1+(t+\frac{1}{2})h)\varepsilon_l^{(s+t)}
\end{eqnarray*}

\textbf{\underline{Brackets involving $HH^{2+6s}(A)$}}:\\

By degree argument, the bracket of $HH^{2+6s}(A)$ with $HH^{2+6t}(A)$ is zero.
\begin{eqnarray*}
 [f_k^{(s)},h_l^{(t)}]&=&\Delta(f_k^{(s)}\cup h_l^{(t)})-\underbrace{\Delta(f_k^{(s)})}_{=0}\cup h_l^{(t)}-f_{k}^{(s)}\cup\Delta(h_{l}^{(t)})\\
 &=&\Delta(\delta_{kl}\psi_0^{(s+t)})-(h+(m-t-1)h)f_k^{(s)}\cup\alpha^{-1}(h_l^{(t)})\\
 &=&\delta_{kl}(h-1+(m-s-t-1)h)\zeta_0-\delta_{kl}(m-t)h\zeta_0\\
 &=&-\delta_{kl}(1+sh)\zeta_0^{(s+t)},\\
 \end{eqnarray*}
 \begin{eqnarray*}
 [f_k^{(s)},\zeta_l^{(t)}]&=&\Delta(f_k^{(s)}\cup\zeta_l^{(t)})-\underbrace{\Delta(f_k^{(s)})}_{=0}\cup\zeta_l^{(t)}-f_k^{(s)}\cup\underbrace{\Delta(\zeta_l^{(t)})}_{=0}\\
 &=&\delta_{l,h-3}(k+1)\Delta(z_{h-3}^{(s+t)})=0,\\
 \end{eqnarray*}

 \begin{eqnarray*}
 [f_k^{(s)},\psi_l^{(t)}]&=&\Delta(f_k^{(s)}\cup\psi_l^{(t)})-\underbrace{\Delta(f_k^{(s)})}_{=0}\cup\psi_l^{(t)}-f_k^{(s)}\cup\Delta(\psi_l^{(t)})\\
 &=&\Delta(f_k^{(s)}\cup\psi_l^{(t)})-(h-1-\frac{l}{2}+(m-t-1)h)f_k^{(s)}\cup\zeta_l^{(t)}\\
&=&\delta_{l,h-3}(k+1)\Delta(\theta_{h-3}^{(s+t+1)})\\&&-\delta_{l,h-3}((m-t)h-1-\frac{h-3}{2})(k+1)z_{h-3}^{(s+t+1)}\\
&=&\delta_{l,h-3}(k+1)(1+\frac{h-3}{2}+(m-s-t-1)h)z_{h-3}^{(s+t+1)}\\
&&-\delta_{l,h-3}((m-t)h-1-\frac{h-3}{2})(k+1)z_{h-3}^{(s+t+1)}\\
&=&-\delta_{l,h-3}(k+1)(1+sh)z_{h-3}^{(s+t+1)},
\end{eqnarray*}

\begin{eqnarray*}
[f_k^{(s)},\varepsilon_l^{(t)}]&=&\Delta(\underbrace{f_k^{(s)})\cup\varepsilon_l^{(t)}}_{=0})-\underbrace{\Delta(f_k^{(s)})}_{=0}\cup\varepsilon_l^{(t)}-\Delta f_k^{(s)}\cup\underbrace{\Delta(\varepsilon_l^{(t)})}_{=0}=0
\end{eqnarray*}

\textbf{\underline{Brackets involving $HH^{3+6s}(A)$}}:\\
We have
\begin{eqnarray*}
 [h_k^{(s)},h_l^{(t)}]&=&\Delta(\underbrace{h_k^{(s)}\cup h_l^{(t)}}_{=0})-\Delta(h_k^{(s)})\cup h_l^{(t)}+h_k^{(s)}\cup\Delta(h_l^{(t)})\\
 &=&-(h+(m-s-1)h)\alpha^{-1}(h_{k}^{(s)})\cup h_l^{(t)}\\&&+(h+(m-t-1)h)h_k^{(s)}\cup\alpha^{-1}(h_l^{(t)})\\
 &=&(s-t)h\alpha^{-1}(h_k^{(s)})\cup h_l^{(t)}=(s-t)h(M_\alpha^{-1})_{kl}\psi_0^{(s+t)},\\
 {}[h_k^{(s)},\zeta_l^{(t)}]&=&\Delta(h_k^{(s)}\cup\zeta_l^{(t)})-\Delta(h_k^{(s)})\cup\zeta_l^{(t)}+\underbrace{h_k^{(s)}\cup\Delta(\zeta_l^{(t)})}_{=0}\\
 {}&=&\delta_{k,\frac{h-3}{2}}\delta_{l,h-3}\Delta(\theta_{h-3}^{(s+t+1)})- (m-s)h\alpha^{-1}(h_k^{(s)})\cup\zeta_l^{(t)}\\
 &=&\delta_{k,\frac{h-3}{2}}\delta_{l,h-3}((1+\frac{h-3}{2}+(m-s-t-1)h)z_{h-3}^{(s+t+1)}\\&&-(m-s)hz_{h-3}^{(s+t+1)})\\
 &=&\delta_{k,\frac{h-3}{2}}\delta_{l,h-3}(-\frac{h+1}{2}-th)z_{h-3}^{(s+t+1)}
  \end{eqnarray*}

We have
\begin{eqnarray*}
[h_k^{(s)},\psi_l^{(t)}]&=&\Delta(\underbrace{h_k^{(s)}\cup\psi_l^{(t)}}_{=0})-\Delta(h_k^{(s)})\cup\psi_l^{(t)}+h_k^{(s)}\cup\Delta(\psi_l^{(t)})\\
&=&-(m-s)h\alpha^{-1}(h_k^{(s)})\cup\psi_l^{(t)}+((m-t)h-1-\frac{l}{2})h_k^{(s)}\cup\zeta_l^{(t)}\\
&=&-(m-s)h\delta_{k,\frac{h-3}{2}}\delta_{l,h-3}\theta_{h-3}^{(s+t+1)}\\&&+((m-t)h-1-\frac{l}{2})\delta_{k,\frac{h-3}{2}}\theta_{h-3}^{(s+t+1)}\\
&=&\delta_{l,h-3}\delta_{k,\frac{h-3}{2}}((s-t)h-\frac{h-1}{2})\theta_{h-3}^{(s+t+1)}
\end{eqnarray*}

\begin{eqnarray*}
[h_k^{(s)},\varepsilon_l^{(t)}]&=&\Delta(\underbrace{h_k^{(s)}\cup\varepsilon_l^{(t)}}_{=0})-\Delta(h_k^{(s)})\cup\varepsilon_l^{(t)}-\Delta h_k^{(s)}\cup\underbrace{\Delta(\varepsilon_l^{(t)})}_{=0}\\
&=&-(m-s)h\alpha^{-1}(h_k^{(s)})\cup\varepsilon_l^{(t)}=0.
\end{eqnarray*}

\textbf{\underline{Brackets involving $HH^{4+6s}(A)$}}:\\

The bracket $[\zeta_k^{(s)},\zeta_l^{(t)}]=\Delta(\zeta_k^{(s)}\cup\zeta_l^{(t)})-\Delta(\zeta_k^{(s)})\cup\zeta_l^{(t)}-\zeta_k^{(s)}\cup\Delta(\zeta_l^{(t)})=0$ because $\Delta$ is zero on $HH^{2+6s}$ and $HH^{4+6s}$.

\begin{eqnarray*}
 [\zeta_{k}^{(s)},\psi_{l}^{(t)}]&=&\Delta(\zeta_{k}^{(s)}\cup\psi_{l}^{(t)})-\underbrace{\Delta(\zeta_{k}^{(s)})}_{=0}\cup\psi_{l}^{(t)}-\zeta_{k}^{(s)}\cup\Delta(\psi_{l}^{(t)}),\\
 &=&\delta_{k,h-3}\delta_{l,h-3}\Delta(\alpha(f_{\frac{h-3}{2}}^{(s+t+1)}))-((m-t)h-1-\frac{l}{2})\zeta_{k}^{(s)}\cup\zeta_{l}^{(t)}\\
 &=&\delta_{k,h-3}\delta_{l,h-3}(m-s-t-1)hf_{\frac{h-3}{2}}^{(s+t+1)}-((m-t)h-1-\frac{h-3}{2})f_{\frac{h-3}{2}}^{(s+t+1)}\\
 &=&\delta_{k,h-3}\delta_{l,h-3}(-sh-\frac{h+1}{2})f_{\frac{h-3}{2}}.
\end{eqnarray*}

\textbf{\underline{The bracket of $HH^{5+6s}(A)$ with $HH^{5+6s}(A)$}}:\\

\begin{eqnarray*}
 [\psi_{k}^{(s)},\psi_{l}^{(t)}]&=&\underbrace{\Delta(\psi_{k}^{(s)}\cup\psi_{l}^{(t)})}_{=0}-\Delta(\psi_{k}^{(s)})\cup\psi_{l}^{(t)}+\psi_{k}^{(s)}\cup\Delta(\psi_{l}^{(t)}),\\
&=&-((m-s)h-1-\frac{k}{2})\zeta_{k}^{(s)}\cup\psi_{l}^{(t)}+((m-t)h-1-\frac{l}{2})\psi_{k}^{(s)}\cup\zeta_{l}^{(t)}\\
&=&\delta_{k,h-3}\delta_{l,h-3}(s-t)h\alpha(f_{\frac{h-3}{2}}^{(s+t+1)}).
\end{eqnarray*}

\newpage

\subsubsection{The Lie derivative $\mathfrak{L}$}
We use the Cartan identity (\ref{Cartan})
to compute the Lie derivative. \\

\textbf{\underline{$HH^{1+6s}(A)$-Lie derivatives}}:\\

From the Cartan identity, we see that
\[
\mathfrak{L}_{\theta_k^{(s)}}=B\iota_{\theta_k^{(s)}}+\iota_{\theta_k^{(s)}}B.
\]
On $\theta_{l,t}$, $\omega_{l,t}$, $\psi_{l,t}$ and $h_{l,t}$, the Connes differential acts by multiplication with $\frac{1}{2}$ degree and taking the preimage under $\iota_{\theta_0}$, and $\iota_{\theta_k^{(s)}}$ acts on them by zero. $B$ acts by zero on $z_l^{(t)}$, $\varepsilon_k^{(t)}$, $\zeta_{l,t}$ and $f_{l,t}$. Since $B$ is degree preserving, this means that $\mathfrak{L}_{\theta_k^{(s)}}$ acts on $\theta_{l,t}$, $\omega_{l,t}$, $\psi_{l,t}$ and $h_{l,t}$ by multiplication with $\frac{1}{2}$ their degree times $z_k^{(s)}$, and on $z_{l,t}$, $\varepsilon_{k,t}$, $\zeta_{l,t}$ and $f_{l,t}$ by multiplication with $z_{k}^{(s)}$ and then multiplication with $\frac{1}{2}$ degree of their product. So we get the following formulas:
\begin{eqnarray*}
 \mathfrak{L}_{\theta_k^{(s)}}(\theta_{l,t})&=&(1+\frac{l}{2}+th)(z_k\theta_l)_{t-s},\\
 \mathfrak{L}_{\theta_k^{(s)}}(z_{l,t})&=&(1+\frac{k+l}{2}+(t-s)h)(z_kz_l)_{t-s},\\
 \mathfrak{L}_{\theta_k^{(s)}}(\omega_{l,t})&=&\delta_{k0}(\frac{1}{2}+t)h\omega_{l,t-s},\\
 \mathfrak{L}_{\theta_k^{(s)}}(\varepsilon_{l,t})&=&\delta_{k0}(\frac{1}{2}+(t-s))h\varepsilon_{l,t-s},\\
 \mathfrak{L}_{\theta_k^{(s)}}(\psi_{l,t})&=&((t+1)h-1-\frac{l}{2})(z_k\psi_l)_{t-s},\\
 \mathfrak{L}_{\theta_k^{(s)}}(\zeta_{l,t})&=&((t-s+1)h-1-\frac{l-k}{2})(z_k\zeta_l)_{t-s},\\
 \mathfrak{L}_{\theta_k^{(s)}}(h_{l,t})&=&\delta_{k0}(t+1)hh_{l,t-s},\\
 \mathfrak{L}_{\theta_k^{(s)}}(f_{l,t})&=&\delta_{k0}(t-s+1)hf_{l,t-s}
\end{eqnarray*}

\newpage

\textbf{\underline{$HH^{2+6s}(A)$-Lie derivatives}}:\\

We compute $\mathfrak{L}_{f_k^{(s)}}$:
\begin{eqnarray*}
 \mathfrak{L}_{f_k^{(s)}}(\theta_{l,t})&=&B(\iota_{f_k^{(s)}}(\theta_{l,t}))-\iota_{f_k^{(s)}}(B(\theta_{l,t}))\\
 &=&B(\delta_{l0}\alpha(f_{k,t-s-1}))-(1+\frac{l}{2}+th))\iota_{f_k^{(s)}}z_{l,t}\\
 &=&\delta_{l0}(t-s)hf_{k,t-s-1}-\delta_{l0}(1+th)f_{k,t-s-1}=-\delta_{l0}(1+sh)f_{k,t-s},\\
 \mathfrak{L}_{f_k^{(s)}}(f_{l,t})&=&B(\underbrace{\iota_{f_k^{(s)}}(f_{l,t})}_{\in HH_{4+6(t-s)}})=0,\\
 \mathfrak{L}_{f_k^{(s)}}(z_{l,t})&=&\delta_{l0}B(f_{k,t-s})=0,\\
 \mathfrak{L}_{f_k^{(s)}}(\omega_{l,t})&=&B(\underbrace{\iota_{f_k^{(s)}}\omega_{l,t}}_{=0})+\iota_{f_k^{(s)}}B(\omega_{l,t})=(\frac{1}{2}+t)h\iota_{f_k^{(s)}}\beta^{-1}(\omega_{l,t})=0,\\
 \mathfrak{L}_{f_k^{(s)}}(\varepsilon_{l,t})&=&B(\underbrace{\iota_{f_k^{(s)}}(\varepsilon_{l,t})}_{=0}),\\
 \mathfrak{L}_{f_k^{(s)}}(\psi_{l,t})&=&B(\iota_{f_k^{(s)}}(\psi_{l,t}))-\iota_{f_k^{(s)}}B(\psi_{l,t})\\
 &=&B(\delta_{l,h-3}(k+1)\theta_{h-3,t-s})-((t+1)h-1-\frac{l}{2})\iota_{f_k^{(s)}}(\zeta_{l,t})\\
 &=&\delta_{l,h-3}(k+1)(1+\frac{h-3}{2}+(t-s)h)z_{h-3,t-s}\\&&-\delta_{l,h-3}(k+1)((t+1)h-1-\frac{l}{2})z_{h-3,t-s}\\
 &=&-\delta_{l,h-3}(k+1)(1+sh)z_{h-3,t-s},\\
 \mathfrak{L}_{f_k^{(s)}}(\zeta_{l,t})&=&B(\iota_{f_k^{(s)}}(\zeta_{l,t}))=B(
k\delta_{l,h-3}z_{h-3,t-s})=0,\\
 \mathfrak{L}_{f_k^{(s)}}(h_{l,t})&=&B(\iota_{f_k^{(s)}}(h_{l,t}))-\iota_{f_k^{(s)}}B(h_{l,t})\\
 &=&B(\delta_{k,l}\psi_{0,t-s})-(t+1)h\iota_{f_k^{(s)}}\alpha^{-1}(h_{l,t})\\
 &=&\delta_{kl}((t-s+1)h-1)\zeta_{0,t-s}-\delta_{kl}(t+1)h\zeta_{0,t-s}\\
 &=&-\delta_{kl}(sh+1)\zeta_{0,t-s}
\end{eqnarray*}

\newpage

\textbf{\underline{$HH^{3+6s}(A)$-Lie derivatives}}:\\

We compute $\mathfrak{L}_{h_k^{(s)}}$:
\begin{eqnarray*}
 \mathfrak{L}_{h_k^{(s)}}(\theta_l^{(t)})&=&B(\underbrace{\iota_{h_k^{(s)}}(\theta_l^{(t)})}_{=0})+\iota_{h_k^{(s)}}B(\theta_{l,t})=(1+\frac{l}{2}+th)\iota_{h_k^{(s)}}z_{l,t}\\
 &=&\delta_{l0}(1+th)h_{k,t-s-1},\\
 \mathfrak{L}_{h_k^{(s)}}(z_{l,t})&=&B(\delta_{l0}h_{k,t-s-1})=\delta_{l0}(t-s)h\alpha^{-1}(h_{k,t-s-1}),\\
 \mathfrak{L}_{h_k^{(s)}}(\omega_{l,t})&=&B\underbrace{\iota_{h_k^{(s)}}(\omega_{l,t})}_{=0}+\underbrace{\iota_{h_k^{(s)}}B(\omega_{l,t})}_{\mbox{cup product in $HH^3(A)\times HH^5(A)$}}=0,\\
 \mathfrak{L}_{h_k^{(s)}}(\varepsilon_{l,t})&=&B\underbrace{\iota_{h_k^{(s)}}\varepsilon_{l,t}}_{=0}=0,\\
 \mathfrak{L}_{h_k^{(s)}}(\psi_{l,t})&=&B(\underbrace{\iota_{h_k^{(s)}}(\psi_{l,t})}_{=0})+h_{k}^{(s)}B(\psi_{l,t})\\
 &=&((t+1)h-1-\frac{l}{2})h_k^{(s)}\zeta_{l,t}=\delta_{k,\frac{h-3}{2}}\delta_{l,h-3}(th+\frac{h+1}{2})\theta_{h-3,t-s},\\
 \mathfrak{L}_{h_k^{(s)}}(\zeta_{l,t})&=&B(\iota_{h_k^{(s)}}(\zeta_{l,t}))=\delta_{k,\frac{h-3}{2}}B(\delta_{l,h-3}\theta_{h-3,t-s})\\&&=\delta_{l,h-3}((t-s)h+\frac{h-1}{2})z_{h-3,t-s},\\
 \mathfrak{L}_{h_k^{(s)}}(h_{l,t})&=&B(\underbrace{\iota_{h_k^{(s)}}(h_{l,t})}_{=0})+\iota_{h_k^{(s)}}B(h_{l,t})=(t+1)h\iota_{h_k^{(s)}}\alpha^{-1}(h_{l,t})\\
 &=&(t+1)h(M_\alpha^{-1})_{lk}\psi_{0,t-s}\\
 \mathfrak{L}_{h_k^{(s)}}(f_{l,t})&=&B(\iota_{h_k^{(s)}}(f_{l,t}))=B(\delta_{kl}\psi_{0,t-s})=\delta_{kl}((t-s+1)h-1)\zeta_{0,t-s}
\end{eqnarray*}

\newpage

\textbf{\underline{$HH^{4+6s}(A)$-Lie derivatives}}:\\

We compute $\mathfrak{L}_{\zeta_k^{(s)}}$:
\begin{eqnarray*}
 \mathfrak{L}_{\zeta_k^{(s)}}(\theta_{l,t})&=&B\iota_{\zeta_k^{(s)}}(\theta_{l,t})-\iota_{\zeta_k^{(s)}}B(\theta_{l,t})=B((z_l\psi_k)_{t-s-1})-\iota_{\zeta_k^{(s)}}(1+\frac{l}{2}+th)z_{l,t}\\
 &=&((t-s)h-1-\frac{k-l}{2})(z_l\zeta_k)_{t-s-1}-(1+\frac{l}{2}+th)(z_l\zeta_k)_{t-s-1}\\
 &=&(-sh-2-\frac{k}{2})(z_l\zeta_k)_{t-s-1},\\
 \mathfrak{L}_{\zeta_k^{(s)}}(\omega_{l,t})&=&0,\\
 \mathfrak{L}_{\zeta_k^{(s)}}(z_{l,t})&=&B\iota_{\zeta_k^{(s)}}(z_{l,t})-\iota_{\zeta_k^{(s)}}\underbrace{B(z_{l,t})}_{=0}\\
 &=&B((z_l\zeta_k)_{t-s-1})=0,\\
 \mathfrak{L}_{\zeta_k^{(s)}}(\psi_{l,t})&=&B\iota_{\zeta_k^{(s)}}(\psi_{l,t})-\iota_{\zeta_k^{(s)}}B(\psi_{l,t})\\
 &=&\delta_{k,h-3}\delta_{l,h-3}B(\alpha(f_{\frac{h-3}{2},t-s-1}))-\iota_{\zeta_k^{(s)}}((t+1)h-1-\frac{l}{2})\zeta_{l,t}\\
 &=&\delta_{k,h-3}\delta_{l,h-3}((t-s)hf_{\frac{h-3}{2},t-s-1}\\&&-((t+1)h-1-\frac{h-3}{2})f_{\frac{h-3}{2},t-s-1})\\
 &=&\delta_{k,h-3}\delta_{l,h-3}(-sh-\frac{h+1}{2})f_{\frac{h-3}{2},t-s-1},\\
 \mathfrak{L}_{\zeta_k^{(s)}}(\varepsilon_{l,t})&=&B\underbrace{\iota_{\zeta_k^{(s)}}(\varepsilon_{l,t})}_{=0}-\iota_{\zeta_k^{(s)}}\underbrace{B(\varepsilon_{l,t})}_{=0}=0,\\
 \mathfrak{L}_{\zeta_k^{(s)}}(\zeta_{l,t})&=&B\iota_{\zeta_k^{(s)}}(\zeta_{l,t})-\iota_{\zeta_k^{(s)}}\underbrace{B(\zeta_{l,t})}_{=0}\\
 &=&\delta_{kh-3}\delta_{l,h-3}B(f_{\frac{h-3}{2},t-s-1})=0,\\
 \mathfrak{L}_{\zeta_k^{(s)}}(h_{l,t})&=&B\iota_{\zeta_k^{(s)}}(h_{l,t})-\iota_{\zeta_k^{(s)}}B(h_{l,t})\\
 &=&
 \delta_{l,\frac{h-3}{2}}\delta_{k,h-3}B(\theta_{h-3,t-s})-(t+1)h\iota_{\zeta_k^{(s)}}\alpha^{-1}(h_{l,t}),\\
 &=&\delta_{l,\frac{h-3}{2}}\delta_{k,h-3}z_{h-3,t-s}((1+\frac{h-3}{2}+(t-s)h)-(t+1)h)\\
 &=&\delta_{l,\frac{h-3}{2}}\delta_{k,h-3}z_{h-3,t-s}(-\frac{h+1}{2}-sh),\\
 \mathfrak{L}_{\zeta_k^{(s)}}(f_{l,t})&=&B\iota_{\zeta_k^{(s)}}(f_{l,t})-\iota_{\zeta_k^{(s)}}\underbrace{B(f_{l,t})}_{=0}\\
 &=&(l+1)\delta_{k,h-3}B(z_{h-3,t-s})=0
\end{eqnarray*}

\newpage 

\textbf{\underline{$HH^{5+6s}(A)$-Lie derivatives}}:\\

We compute $\mathfrak{L}_{\varepsilon_k^{(s)}}$:
\begin{eqnarray*}
\mathfrak{L}_{\varepsilon_k^{(s)}}(\theta_{l,t})&=&B(\iota_{\varepsilon_k^{(s)}}(\theta_{l,t}))+\iota_{\varepsilon_k^{(s)}}B(\theta_{l,t})\\
&=&B(-\delta_{l0}\beta(\varepsilon_{k,t-s-1}))+(1+\frac{l}{2}+th)\iota_{\varepsilon_k^{(s)}}(z_{l,t})\\
&=&-\delta_{l0}(\frac{1}{2}+t-s-1)h\varepsilon_{k,t-s-1}+(1+th)\delta_{l0}\varepsilon_{k,t-s-1}\\
&=&((s+\frac{1}{2})h+1)\delta_{l0}\varepsilon_{k,t-s-1},\\
\mathfrak{L}_{\varepsilon_k^{(s)}}(z_{l,t})&=&B(\iota_{\varepsilon_k^{(s)}}(z_{l,t}))=B(\varepsilon_{k,t-s-1})=0,\\
\mathfrak{L}_{\varepsilon_k^{(s)}}(\omega_{l,t})&=&B(\iota_{\varepsilon_k^{(s)}}(\omega_{l,t}))+\iota_{\varepsilon_k^{(s)}}B(\omega_{l,t})\\
&=&B(\delta_{kl}\psi_{0,t-s-1})+(\frac{1}{2}+t)h\iota_{\varepsilon_k^{(s)}}\beta^{-1}(\omega_{l,t})\\
&=&\delta_{kl}((t-s)h-1)\zeta_{0,t-s-1}-\delta_{kl}(\frac{1}{2}+t)h\zeta_{0,t-s-1}\\&=&-\delta_{kl}(1+(\frac{1}{2}+s)h)\zeta_{0,t-s-1},\\
\mathfrak{L}_{\varepsilon_k^{(s)}}(\psi_{l,t})&=&B(\underbrace{\iota_{\varepsilon_k^{(s)}}(\psi_{l,t})}_{=0})+\iota_{\varepsilon_k^{(s)}}B(\psi_{l,t})\\
&=&((t+1)h-1-\frac{l}{2})\iota_{\varepsilon_k^{(s)}}\zeta_{l,t}=0,\\
\mathfrak{L}_{\varepsilon_k^{(s)}}(\varepsilon_{l,t})&=&B(\iota_{\varepsilon_k^{(s)}}(\varepsilon_{l,t}))=B(-(M_\beta)_{kl}\zeta_{0,t-s-1})=0,\\
\mathfrak{L}_{\varepsilon_k^{(s)}}(\zeta_{l,t})&=&B(\underbrace{\iota_{\varepsilon_k^{(s)}}\zeta_{l,t}}_{=0})=0,\\
\mathfrak{L}_{\varepsilon_k^{(s)}}(h_{l,t})&=&B(\iota_{\varepsilon_k^{(s)}}(h_{l,t}))+(t+1)h\iota_{\varepsilon_k^{(s)}}\alpha^{-1}(h_{l,t})=0,\\
\mathfrak{L}_{\varepsilon_k^{(s)}}(f_{l,t})&=&B(\underbrace{\iota_{\varepsilon_k^{(s)}}(f_{l,t})}_{=0})=0
\end{eqnarray*}

\newpage

We compute $\mathfrak{L}_{\psi_k^{(s)}}$:

\begin{eqnarray*}
 \mathfrak{L}_{\psi_k^{(s)}}(\theta_{l,t})&=&B\underbrace{\iota_{\psi_k^{(s)}}(\theta_{l,t})}_{=0}+\iota_{\psi_k^{(s)}}B(\theta_{l,t})\\
 &=&\iota_{\psi_k^{(s)}}z_{l,t}(1+\frac{l}{2}+th)=(z_l\psi_k)_{t-s-1}(1+\frac{l}{2}+th),\\
 \mathfrak{L}_{\psi_k^{(s)}}(z_{l,t})&=&B\iota_{\psi_k^{(s)}}(z_{l,t})+\iota_{\psi_k^{(s)}}\underbrace{B(z_{l,t})}_{=0}\\
 &=&B((z_l\psi_k)_{t-s-1})=((t-s)h-1-\frac{k-l}{2})(z_l\zeta_k)_{t-s-1},\\
 \mathfrak{L}_{\psi_k^{(s)}}(\omega_{l,t})&=&B\underbrace{\iota_{\psi_k^{(s)}}(\omega_{l,t})}_{=0}+\iota_{\psi_k^{(s)}}B(\omega_{l,t})\\
 &=&(\frac{1}{2}+t)h\iota_{\psi_k^{(s)}}\beta^{-1}(\omega_{l,t})=0,\\
 \mathfrak{L}_{\psi_k^{(s)}}(\psi_{l,t})&=&B\underbrace{\iota_{\psi_k^{(s)}}(\psi_{l,t})}_{=0}+\iota_{\psi_k^{(s)}}B(\psi_{l,t})\\
 &=&((t+1)h-1-\frac{l}{2})\iota_{\psi_k^{(s)}}\zeta_{l,t}\\
 &=&\delta_{k,h-3}\delta_{l,h-3}(\underbrace{(t+1)h-1-\frac{h-3}{2}}_{=th+\frac{h+1}{2}})\alpha(f_{\frac{h-3}{2},t-s-1}),\\
 \mathfrak{L}_{\psi_k^{(s)}}(\varepsilon_{l,t})&=&B\underbrace{\iota_{\psi_k^{(s)}}(\varepsilon_{l,t})}_{=0}+\iota_{\psi_k^{(s)}}\underbrace{B(\varepsilon_{l,t})}_{=0}=0\\
 \mathfrak{L}_{\psi_k^{(s)}}(\zeta_{l,t})&=&B\iota_{\psi_k^{(s)}}(\zeta_{l,t})+\iota_{\psi_k^{(s)}}\underbrace{B(\zeta_{l,t})}_{=0}\\
 &=&\delta_{k,h-3}\delta_{l,h-3}B(\alpha(f_{\frac{h-3}{2},t-s-1}))\\
 &=&\delta_{k,h-3}\delta_{l,h-3}(t-s)hf_{\frac{h-3}{2},t-s-1},\\
 \mathfrak{L}_{\psi_k^{(s)}}(h_{l,t})&=&B\underbrace{\iota_{\psi_k^{(s)}}(h_{l,t})}_{=0}+\iota_{\psi_k^{(s)}}B(h_{l,t})=\iota_{\psi_k^{(s)}}\alpha^{-1}(h_{l,t})(t+1)h\\
 &=&\delta_{k,h-3}\delta_{l,\frac{h-3}{2}}(t+1)h\theta_{h-3,t-s},\\
 \mathfrak{L}_{\psi_k^{(s)}}(f_{l,t})&=&B_{\psi_k^{(s)}}(f_{l,t})+\iota_{\psi_k^{(s)}}\underbrace{B(f_{l,t})}_{=0}\\
 &=&(l+1)\delta_{k,h-3}B(\theta_{h-3,t-s})\\&=&(l+1)(1+(t-s)h+\frac{h-3}{2})\delta_{k,h-3}z_{h-3,t-s}
\end{eqnarray*}

\newpage

\textbf{\underline{$HH^{6+6s}(A)$-Lie derivatives}}:\\

$B$ acts on $\theta_{l,t},\,\omega_{l,t},\,\psi_{l,t}\,$ and $h_{l,t}$ by multiplication with $\frac{1}{2}$ degree and taking the preimage under $\iota_{\theta_0}$. On $z_{l,t},\,\varepsilon_{l,t},\,\zeta_{l,t}$ and $f_{l,t}$, $B$ acts by zero. Since the spaces $U$, $U^*$, $K$, $K^*$, $Y$ and $Y^*$ are $z_k$-invariant and $z_k^{(s)}$ has degree $k-2sh$, $\mathfrak{L}_{z_k^{(s)}}$ acts on  $\theta_{l,t},\,\omega_{l,t},\,\psi_{l,t}\,$ and $h_{l,t}$ by multiplication with $\frac{k}{2}-sh$ and taking the preimage under $\iota_{\theta_0}$ and multiplication with $z_k^{(s)}$, and on $z_{l,t},\,\varepsilon_{l,t},\,\zeta_{l,t}$ and $f_{l,t}$ it acts by zero. We have the following formulas:

\begin{eqnarray*}
 \mathfrak{L}_{z_k^{(s)}}(\theta_{l,t})&=&(\frac{k}{2}-sh)(z_k\theta_l)_{t-s},\\
 \mathfrak{L}_{z_k^{(s)}}(z_{l,t})&=&0,\\
 \mathfrak{L}_{z_k^{(s)}}(\omega_{l,t})&=&-\delta_{k0}sh\beta^{-1}(\omega_{l,t-s}),\\
 \mathfrak{L}_{z_k^{(s)}}(\psi_{l,t})&=&(\frac{k}{2}-sh)(z_k\zeta_l)_{t-s},\\
 \mathfrak{L}_{z_k^{(s)}}(\varepsilon_{l,t})&=&0,\\
 \mathfrak{L}_{z_k^{(s)}}(\zeta_{l,t})&=&0,\\
 \mathfrak{L}_{z_k^{(s)}}(h_{l,t})&=&(\frac{k}{2}-sh)\alpha^{-1}(h_{l,t-s}),\\
 \mathfrak{L}_{z_k^{(s)}}(f_{l,t})&=&0
\end{eqnarray*}

Now we compute $\mathfrak{L}_{\omega_k^{(s)}}$:

We observe that $\iota_{\omega_k^{(s)}}(\varepsilon_{l,t})=\delta_{kl}\psi_{0,t-s}$,  $\iota_{\omega_k^{(s)}}(z_{l,t})=\delta_{l0}\omega_{k,t-s}$, and
\[
 \iota_{\omega_k^{(s)}}(\theta_{l,t})=\iota_{\omega_k^{(s)}}(\omega_{l,t})=\iota_{\omega_k^{(s)}}(\psi_{l,t})=\iota_{\omega_k^{(s)}}(\zeta_{l,t})=\iota_{\omega_k^{(s)}}(h_{l,t})=\iota_{\omega_k^{(s)}}(f_{l,t})=0.
\]
Then we have
\begin{eqnarray*}
 \mathfrak{L}_{\omega_k^{(s)}}(\varepsilon_{l,t})&=&B\iota_{\omega_k^{(s)}}(\varepsilon_{l,t})=\delta_{kl}B(\psi_{0,t-s})\\
 &=&\delta_{kl}((t-s+1)h-1)\zeta_{0,t-s},\\
 \mathfrak{L}_{\omega_k^{(s)}}(z_{l,t})&=&B\iota_{\omega_k^{(s)}}(z_{l,t})=\delta_{l0}B(\omega_{k,t-s})\\
 &=&\delta_{l0}(\frac{1}{2}+t-s)h\beta^{-1}(\omega_{k,t-s}),\\
 \mathfrak{L}_{\omega_k^{(s)}}(\theta_{l,t})&=&\iota_{\omega_k^{(s)}}B(\theta_{l,t})=(1+\frac{l}{2}+th)\iota_{\omega_k^{(s)}}z_{l,t}\\
 &=&\delta_{l0}(1+th)\omega_{k,t-s},
\end{eqnarray*}
and 
\[
 \mathfrak{L}_{\omega_k^{(s)}}(\omega_{l,t})= \mathfrak{L}_{\omega_k^{(s)}}(\psi_{l,t})=\mathfrak{L}_{\omega_k^{(s)}}(\zeta_{l,t})= \mathfrak{L}_{\omega_k^{(s)}}(h_{l,t})= \mathfrak{L}_{\omega_k^{(s)}}(f_{l,t})=0.
\]

\end{document}